# Temporal truth and bivalence: an anachronistic formal approach to Aristotle's *De Interpretatione* 9

Luiz Henrique Lopes dos Santos

Universidade de São Paulo

*Abstract*

Regarding the famous Sea Battle Argument, which Aristotle presents in *De Interpretatione* 9, there has never been a general agreement not only about its correctness but also, and mainly, about what the argument really is. According to the most natural reading of the chapter, the argument appeals to a temporal concept of truth and concludes that not every statement is always either true or false. However, many of Aristotle's followers and commentators have not adopted this reading. I believe that it has faced so much resistance for reasons of hermeneutic charity: denying the law of universal bivalence seems to be overly disruptive to logical orthodoxy – the kind of logical orthodoxy represented by what we now call classical propositional logic, much of which Aristotle clearly supports in many texts. I intend to show that the logical-semantic theses that the traditional reading finds in *De Interpretatione* 9 are much more conservative than they may seem to be at first glance. First, I will show that they complement, and do not contradict in any way, the orthodox definitions of the concepts of truth and statement that Aristotle advances in other texts. Second, by resorting in an anachronistic vein to concepts and methods peculiar to contemporary logic, I will show that a trivalent modal semantics conforming to those theses can be built for a standard formal language of the classical propositional calculus. It is remarkable that reasonable concepts of logical truth and logical consequence that may be defined on the basis of this trivalent modal semantics are coextensive with their orthodox counterparts, the concepts of tautology and tautological consequence of classical bivalent and extensional semantics.

*Keywords*: Aristotle, *De Interpretatione* 9, Sea Battle Argument, law of bivalence, temporal truth, trivalent semantics.

*Introduction*

The so-called Sea Battle Argument, which Aristotle presents in *De Interpretatione* 9, was one of the most famous topics of discussions about determinism in Antiquity and the Middle Ages. As a matter of fact, it was as famous as it was and has always been highly controversial. There has never been a general agreement not only about its formal validity and the truth of its premises, but also, and mainly, about what the argument itself really is, that is, about what its premises really are and what its conclusion really is.

This may seem curious to a lay reader of *De Interpretatione* 9, for the letter of the chapter naturally and strongly suggests what should be its correct exegesis. The most natural reading of the chapter is the one now improperly qualified as *traditional*. There it finds the exposition of an argument by reduction to absurdity of what we now call the *law of universal bivalence*: every statement is either true or false. According to the traditional reading, Aristotle takes this argument to be an impeccable foundation of the following contention: bivalence is a universal attribute of statements about the present and the past, but not of statements about the future – for it is not an attribute of statements that affirm or deny the reality of contingent future facts.

I take this reading to be *improperly* called the traditional one because there is good historical evidence that it has not been adopted by most of Aristotle's ancient and medieval followers and commentators.[1] It has also been contested by many highly-regarded contemporary commentators of *De Interpretatione* 9.[2]

I believe that the traditional reading faces so much resistance because denying the law of universal bivalence seems to be overly disruptive to logical orthodoxy – the kind of logical orthodoxy represented by what we now call classical propositional logic, much of which Aristotle clearly supports in many texts. In particular, denying the universal validity of the law of bivalence would be incompatible with the definition of truth that

---

[1] In the 6th century AD, Simplicius testified that this reading was rejected by those he called the Peripatetics (cf. Simplicius, *Commentary to Aristotle's Categories*, 406, 5-16; 406, 34-407, 14). Around the same time, Boethius, who would become one of Aristotle's most influential commentators in the Middle Ages, accused the Stoics of wrongly ascribing to Aristotle the thesis that statements about contingent futures are neither true nor false (cf. Boethius, *Second Commentary to Aristotle's De Interpretatione*, 208, 1-11).

[2] For example, Anscombe 1956, Strang 1960, Rescher 1963, Fine 1984, and Judson 1988.

Aristotle proposes in the *Metaphysics*;[3] it would compromise the universal validity of the law of excluded middle, also expressly stated in the *Metaphysics*[4]; and it would compromise the essential link between the concept of statement and the attribute of being true or false, a link that Aristotle expressly mentions in defining this concept in *De Interpretatione* 4.[5]

In a recent long and detailed article, I sustained the traditional reading against its competitors.[6] However, my concern here is not the issue of what would be the correct reading of the chapter, but rather the issue of whether the logical-semantic theses that this reading finds in *De Interpretatione* 9 are really so weird as to deserve so much resistance. It is remarkable that even some contemporary logicians who adopted the traditional reading of the chapter – such as Quine, for example – believe there is no way to systematize the logical-semantic theses they find stated or implicated in the Sea Battle Argument in such a way as to build on them a reasonable logical semantics for propositional logic.[7]

Contrary to these assessments, I intend to show that *De Interpretatione* 9, according to the traditional reading, defends a set of logical-semantic theses that is much less disruptive than it may seem at first sight - either to logical orthodoxy or to the whole body of Aristotelian texts.

*The Sea Battle Argument: premises, conclusions, and implications*

Let us begin by sketching broadly the argumentative movement that the traditional reading finds in *De Interpretatione* 9, in order to identify the logical-semantic theses that the Sea Battle Argument presupposes, intends to substantiate, or implies. According to this reading, in the first stage of the argument (18a34-b25), Aristotle aims to prove that the law of universal bivalence implies a thesis that, in the second stage (18b26-19a22), he aims to show to be absurd. From these partial conclusions, he infers, by *modus tollens*, the negation of the law of universal bivalence.

---

[3] Cf. *Metaphysics* IV 7, 1011b25-27; IX 10, 1051b3-5.
[4] Cf. *Metaphysics* IV 7.
[5] Cf. *De Interpretatione* 4, 17a1-5; also, *Categories* 4, 2a4-10.
[6] Cf. Santos 2021.
[7] Quine qualifies as "a fantasy" Aristotle's acceptance of the law of excluded middle while denying the law of universal bivalence (cf. Quine 1953, p. 65).

The thesis that, in the first stage of the argument, Aristotle intends to show to be implied by the principle of bivalence is determinism: everything that exists or happens has always been determined to exist or happen, exactly when it exists or happens, for whatever reasons (logical, physical, metaphysical, or of any other kind). Conforming to the broader concept of necessity, usually adopted by Aristotle, determinism is the thesis that everything that exists and happens exists or happens by necessity.

In the second stage of the argument, Aristotle does not demonstrate, in the proper sense of the word, the falsity of determinism, but intends to show that it has unacceptable consequences, in light of something that he believes our experience of the world reveals to be incontestable. For him, it is self-evident that there are facts that, although it is possible for them to occur in the future, will never occur, facts that are now neither determined to occur nor determined not to occur in the future, facts whose future reality or unreality is contingent – for instance, possible deliberate human actions. Now, since determinism is false and follows from the law of universal bivalence, Aristotle concludes that this law does not hold for all statements, for it does not hold for those that affirm the contingent reality of future facts.

In the last stage of the argument (19a23-36), to point up *why*, and not only *that*, the law of bivalence is not universally valid, Aristotle turns his reductive argument into an equivalent one in the form of a *modus ponens*. The way he justifies the first premise of this argument, that is, the way he explains why statements about contingent futures are neither true nor false commits him to hold that only *temporal* definitions of truth and falsity are genuinely germane to *temporal statements*, i.e., statements that assert the obtaining of state of affairs at definite times. I advocate that these temporal definitions specify and complement, and do not contradict in any way, the *temporally neutral* definitions proposed in the *Metaphysics*.

Finally, Aristotle lays down a law whose significance in the context of the Sea Battle Argument I think has not yet been sufficiently stressed (19a36-39). He replaces the law of universal bivalence with a weaker analogue, which may be called the *law of universal weak bivalence*: although not every statement is at any moment true or false, at any moment *every* statement necessarily is *or will be* true or false.

Here it is not relevant to scrutinize how Aristotle argues that determinism is unacceptable, but rather to look into the first and last stages of the Sea Battle Argument. Foremost, it should be noted that the alethic modalities involved in *De Interpretatione* 9, according to the traditional reading, are *temporally relative* modalities, of which Aristotle also makes use in other texts, such as *Nicomachean Ethics*, *Rhetoric*, and *On the Heavens*.[8]

In the temporally relative sense of modal terms, something is said to be necessary, possible or impossible *in relation to moments of time*. In this sense, it is necessary at a moment m that a fact is real at a moment m* if and only if it is determined at m that this fact is real at m*; it is impossible at a moment m that a fact is real at a moment m* if and only if it is determined at m that the fact is not real at m*; temporally relative possibility and temporally relative contingency are defined similarly, in compliance with the usual cross-definitions of alethic modalities.

Temporally relative modalities are theoretically fruitful to the extent that at least in principle there can be facts that at a given moment are determined to be real at a given time, but at another moment were not yet determined to be real at that time. For instance, from an indeterminist point of view, two years ago it was not necessary for me to be out of São Paulo now; however, two minutes ago, when I was more than a thousand miles far from São Paulo, it has certainly become necessary for me to be out of São Paulo now.

In conjunction with the principle of non-contradiction, the definitions of temporally relative modalities imply the so-called law of the necessity of the present: if a fact is real at a definite moment, then it is necessary at this moment that it occurs at this very moment. By the principle of non-contradiction, if something is real at a time, then it is impossible for it not to be real at that same time.[9]

Assuming the past to be irreversible, the definitions of temporally relative modalities also imply the so-called law of the necessity of the past: if a fact was real at a moment m before moment m*, then it is necessary at m* that it was real at m. From the moment a state of affairs obtains, it becomes necessary that it has obtained at that moment. From the time when the Greeks sacked Troy, for example, it became forever impossible for the Greeks not to have sacked Troy at that time.[10]

---

[8] Cf. *Nicomachean Ethics* VI 2, 1139b5-11; *Rhetoric* III 17, 1417b38-1418a5; *On the Heavens* I 12 *passim*.

[9] According to the traditional reading, Aristotle states the law of the necessity of the present in *De Interpretatione* 9, 19a23-27. For the sake of argumentative fair play, it should be noted that the sense of this passage is as controversial as the traditional reading itself.

The Sea Battle Argument has explicit and implicit premises. The explicit ones are the principle of non-contradiction, the definition of truth, and indeterminism. The implicit ones are the modal laws of the necessity of the present and the past, and another modal law, actually a very trivial one: what follows from something necessary is also necessary.[11]

Let us now sketch in a free manner the Aristotelian proof that the law of universal bivalence implies determinism. Suppose that every statement is now true or false, and consider the statement that a fact will occur at some future time –for example, the statement that there will be a sea battle tomorrow. If that statement is now true, then, by the necessity of the present, it is now necessary for it to be now true. But if it is now true, then there will be a sea battle tomorrow, by the definition of truth. Since what necessarily follows from what is necessary is also necessary, if the statement that there will be a sea battle tomorrow is now necessarily true, then it is now necessary that there will be a sea battle tomorrow, it is now already determined that there will be a sea battle tomorrow.

We can prove similarly that if this statement is now false, then it is now necessary that there will be no sea battle tomorrow. Hence, if the statement that there will be a sea battle tomorrow is now either true or false, it follows that either it is now necessary that there will be a sea battle tomorrow, or it is now necessary that there will be no sea battle tomorrow.

There is no doubt that this argument is formally valid and can be made completely general, holding for any moment of time, any statement, and any occurrence of any fact at any moment. It is reasonable to assume that any occurrence of any fact at any moment can at least in principle be stated in some possible language. Therefore, if every statement is true or false, then everything that happens at a moment m has always been determined to happen at m, and everything that does not happen at a moment m has always been determined not to happen at m. The assumption of the universal validity of the law of bivalence implies that nothing happens or fails to happen contingently. The consequent being plainly false, Aristotle goes on, the initial assumption must be plainly false. Not every statement is always true or false.

---

[10] The law of the necessity of the past is clearly implicated in *Nicomachean Ethics* VI 2, 1139b5-11.
[11] This law is stated in *Eudemian Ethics* II 6, 1223a1.

From an indeterminist point of view, the first and second stages of the Sea Battle Argument are enough to justify the denial of the law of universal bivalence. Nonetheless, they do not make clear *why* the statement of a future contingent occurrence of a state of affairs must be neither true nor false. Paradoxically, the answer to this question can be extracted from the definition of truth that Aristotle lays down in the *Metaphysics*, paradoxically the same definition that seems to imply the law of universal bivalence.

A literal and unbiased reading of the definitions of truth and falsity formulated in the *Metaphysics* reveals that they are indeed temporally neutral. They assert that a statement *is* true if and only if what it states to be real *is* actually real, and it *is* false if and only if what it states to be real *is* actually unreal. All we can conclusively infer from this definition is the logical equivalence between the truth of a statement and the reality of what it says to be real, as well as the logical equivalence between the falsity of a statement and the unreality of what it says to be real.

At the same time, Aristotle remarks that while the logical relation between the truth or falsity of a statement and the reality or unreality of what it affirms to be real is obviously symmetrical, there is between them an *asymmetrical causal relation*. As he points out, it is not because it is true to say that you are pale that you are pale, but it is because you are pale that it is true to say that you are pale.[12] Reality is *the* cause of truth, unreality is *the* cause of falsity, not vice-versa.

In the context of *De Interpretatione* 9, from the conjunction of logical equivalence and causal asymmetry between truth and reality, it follows that statements about contingent futures can be neither true nor false. In fact, if it is not determined at a given moment whether a necessary cause will be real or unreal, then it is not determined at that moment whether its effect will be real or unreal. Therefore, if now it is not determined that a sea battle will take place tomorrow, then now it is not determined that the statement that a sea battle will take place tomorrow is now true; hence, by the necessity of the present, now the statement is not true. Similarly, if now it is not determined that a sea battle will not take place tomorrow, then now it is not determined that the statement that a sea battle will take place tomorrow is false; therefore, by the necessity of the present, now the statement is not false. Statements about indeterminate futures are now neither true nor false.

---

[12] Cf. *Metafísica* IX 10, 1051b6-9

So, the conjunction of logical equivalence and causal asymmetry between truth and reality imposes the adoption of the following *modal* definition of temporal truth: a statement is true at a given moment if and only if at that moment *it is necessary* to be real what the statement says to be real. And similarly concerning falsity.

Laying these definitions of temporal truth and falsity is tantamount to avowing that the so-called Tarski's T-scheme, and its analogue concerning falsity, are not universally valid in the domain of temporal statements. In this domain, not every statement S is such that, for any moment, S is true (false) at that moment if and only if S (not S). The T-scheme and its analogue concerning falsity hold for statements about the present, the past, and the necessary future, but not for statements about contingent futures. So, Aristotle is free to refuse the equivalence between the law of excluded middle and the law of universal bivalence, an equivalence that is trivially implied by the T-scheme and its analogue concerning falsity. In other words, he is free to refuse universal bivalence without prejudice to the universality of the law of excluded middle.

After explaining why statements about contingent futures are neither true nor false, Aristotle postulates a weaker version of the law of universal bivalence. The *law of universal weak bivalence* asserts that every statement either is *or will be* true or false. Statements about contingent futures are neither true nor false, but they will necessarily be true or false at the right moment. The statement that there will be a sea battle tomorrow is now neither true nor false, but it will certainly be true or false tomorrow at midnight at the latest.

Assuming that the Sea Battle Argument is a good one, the definition of an appropriate semantics for propositional logic must satisfy a number of conditions. We saw that a remarkable one is that the axes of this semantics must be temporally relative concepts of truth and falsity, concepts that are to be temporally relative in the strictest sense: at least in principle, statements can be true (false) at a given moment and not be true (false) at another moment. For the importance and novelty of this condition to be properly measured, it is necessary to clarify how it is to be precisely understood. To do so, it is necessary to dispel one among the many ambiguities that pervade the use of the Greek expression *logos apophantikos*, which I translate by the word "statement", ambiguities that also pervade the use of this English word.

I call a *sentence* every symbol capable of conveying statements. We can distinguish two kinds of sentences capable of conveying affirmations and negations that states of affairs obtain at definite times. On the one hand, there are sentences in which the reference to definite times is made using deictic symbolic resources – such as verbal tenses and adverbial expressions like "now", "tomorrow", and "in the future". I label these sentences, like "Socrates is sitting" and "There will be a sea battle tomorrow", as *temporally relative sentences*. Utterances of the same temporally relative sentence at different moments can, without changing the meanings of its parts, have different truth conditions, since these utterances affirm or deny the occurrence of the same state of affairs, but possibly at different times.

On the other hand, there are sentences in which the reference to definite times is made using expressions that always refer to the same times at all moments when they are uttered, such as expressions for dates and those referring to particular events. I label these sentences, like "There will be a sea battle on September 14, 2022" and "Socrates' death precedes Plato's death", as *temporally absolute sentences*. For simplicity, I will consider only sentences in which no non-temporal deixes occur. With that proviso, all utterances of a temporally absolute sentence have the same truth condition, since all of them affirm or deny that the same state of affairs occurs at the same definite time.

The ambiguity of the word "statement" that matters here accounts for the fact that there are two acceptable but incompatible answers to the question of which statement the utterance of a temporally relative sentence conveys. In a sense, we can say that the sentence "There will be a sea battle tomorrow" has the same meaning at all moments when it is uttered. This common meaning may be called the statement it expresses, namely, the statement that there will be a sea battle the day after the moment when the sentence is uttered, whatever that moment may be. It is in this sense of the word "statement", *the temporally relative sense*, that Aristotle uses the word in *Categories* 5, for instance.

Nevertheless, it is also acceptable to say that a statement is defined by its *whole* truth conditions. In this *temporally absolute sense* of the word "statement", different utterances of the sentence "There will be a sea battle tomorrow" can convey different statements. When that sentence is uttered on September 13, 2022, it conveys the same statement as the sentence "There will be a sea battle on September 14, 2022"; when it is uttered the day after September 13, 2022, it states the same as the sentence "There will be a sea battle on September 15, 2022". I call a *temporally relative statement* what a temporally relative sentence conveys in the temporally relative sense; I call a *temporally absolute statement* what any sentence conveys in the temporally absolute sense.

It is nothing but a triviality to admit that temporally relative statements can have different truth values at different moments. However, it is not a triviality what Aristotle's refusal of the law of strong bivalence requires. It requires that also the logical-semantic treatment for temporally absolute statements resorts to temporally relative concepts of truth and falsity. Before the Greeks sacked Troy, the temporally absolute statement conveyed by an utterance now of the sentence "The Greeks sacked Troy" was neither true nor false, on the plausible assumption that, before it occurred, the sacking of Troy could have been prevented by contingent choices of the Greek commanders or some chance events. From the moment of the sacking, however, the necessity of the present and the past ensured that this same statement became forever true.

Now, supporters of the traditional reading like myself must assume that it is in the temporally absolute sense that Aristotle uses the word "statement" when formulating the law of weak bivalence in *De Interpretatione* 9. It is not difficult to see that, from an indeterminist point of view, the law of weak bivalence may not hold for some temporally relative statements.

From this point of view, for example, it is plausible to assume that as long as there are human beings, the occurrence or non-occurrence of a sea battle, which depend on acts of choice and many kinds of chance events, will never be determined to happen or not to happen a year before it happens or not happen. Insofar as it is at least possible that, as Aristotle believes, the human species will exist throughout infinite time, it follows that at no moment the temporally relative statement conveyed by the sentence "There will be a sea battle within a year" will be true or false, since at no moment it will be determined whether a sea battle will or will not occur within a year. As a consequence, for the sake of hermeneutic charity, the traditional reading of *De Interpretatione* 9 is committed to postulating that statements at stake in the context of the chapter are temporally absolute statements.

However, from the indeterminist point of view, it seems that the law of weak bivalence does not necessarily hold even for all temporally absolute statements. Assuming that time is infinite, as Aristotle does, it seems not to hold for temporally absolute statements that assert the reality at some indefinite moment in the future of a state of affairs that will never be real in the future, although it will be *always* possible for it to be real in the future. The future reality of this state of affairs will be at every moment contingent, so that the statement that it will be real in the future will remain forever neither true nor false.

Therefore, the law of weak bivalence implies that no state of affairs can be such that its future occurrence remains eternally contingent without ever becoming either real or impossible. As a matter of fact, it is equivalent to a weak version of the so-called principle of plenitude.

The strongest version of the principle of plenitude asserts that whatever is possible at a moment is either real at that moment or will eventually become real in the future. In *De Interpretatione* 9, Aristotle expressly rejects this version (19a7-18). It is obvious, he says, that many states of affairs are now possible to be real in the future, but will never become real. This cloak, he says, now may be and may not be cut in the future, but it may wear out and cease to exist before being cut.

The weak version of the principle of plenitude is a little less ambitious. It does not imply that something which is *now* merely possible must eventually become real, but it claims that nothing can remain *eternally* possible without ever becoming real. Excluded by the strongest version of the principle of plenitude, the possibility of being cut the Aristotelian cloak that will never be actually cut, for instance, is rescued by its weak version and so does not challenge the law of weak bivalence. It is now possible for the cloak to be cut; however, from the moment it wears out and ceases to exist before being cut, it will become impossible for it to be cut, and from that moment the statement that the cloak will be cut will be always false. There are good textual reasons to believe that Aristotle admits some weak version of the principle of plenitude, and there are also good reasons to believe that he admits precisely the weak version that is equivalent to the law of weak bivalence.[13]

*DI9 Formal Semantics*

In sum, according to the traditional reading of *De Interpretatione* 9, the Sea Battle Argument presupposes or implies the following logical-semantic theses:

(1) truth and falsity belong primarily to temporally absolute statements;

(2) truth and falsity are temporally relative attributes of temporal statements;

---

[13] Cf. Santos 2021, 113-125.

(3) at each moment, there are for a temporal statement, in principle, three alternatives: to be true, to be false, or to be neither true nor false;

(4) a statement is true (false) at a definite moment if and only if it is necessary for it at that moment to be true (false) at some moment;

(5) statements that are true (false) at a definite moment will remain true (false) forever after that moment;

(6) at every moment every statement is *or will be* true or false.

A deliberate and open anachronistic logical exercise can reveal how close these theses are to classical propositional logic. Indeed, they can underpin the definition of a reasonable formal semantics (which I will call *DI9 semantics*) for a standard formal language of what we now call classical propositional calculus; and it can be proven that reasonable concepts of logical truth and logical consequence defined on the basis of DI9 semantics are coextensive with their classical counterparts, the concepts of tautology and tautological consequence of classical semantics.

Theses (3) and (4) above imply that DI9 semantics should be a *trivalent modal* semantics. It will be defined through a possible worlds strategy. It must entail a definition of truth that makes true at a given moment in a given world exactly those statements that are true or will be true in all possible worlds compatible with the totality of what is real or has already been real at that moment in this world, that is, in all possible continuations of the totality of real occurrences of states of affairs up to that moment in this world.[14]

Let L be a standard language of the classical propositional logic whose primitive connectives are those of negation (~) and disjunction (∨). The set of formulas of L is recursively defined in the usual way from an infinite number of atomic formulas and those connectives.

---

[14] The idea of formally dealing with temporal truth and truth value gaps of statements about the future by resorting to a possible worlds strategy, similar in some important respects to the one employed here in defining DI9 semantics, goes back to Thomason 1970. However, Thomason and his followers do not deal with standard languages of the classical propositional logic, whose formulas are intended to stand for temporally absolute statements. They are concerned with tense languages, whose formulas are intended to stand for temporally relative statements and may include tense operators, corresponding to ordinary expressions like "it was the case that", "it will always be the case that", etc.

A DI9 interpretation for L should be an assignment of one of the values T, F and O to each ordered pair (A, j), A being a formula of L and j being a real number. The ordered set of real numbers will be used as a formal proxy for the ordered set of moments of time (in deference to Aristotle's belief in the continuity and infinity of time). Values T and F will be called *truth values*. It is worth noting that the value 0 is not meant to be a third *truth* value. It is meant to represent formally the (provisory) *absence* of truth value.

The initial step in defining the concept of DI9 interpretation is to define the concepts of DI9 valuation for L and DI9 classical interpretation for L.

*Definition* 1. A *DI9 valuation* for L is a function $\alpha$ from the Cartesian product of the set of atomic formulas of L and the set of real numbers to the set of values $\{T, F, 0\}$ such that, for any atomic formula A of L,

1) for all numbers j and h such that $j < h$, if $\alpha(A, j) = T$, then $\alpha(A, h) = T$; if $\alpha(A, j) = F$, then $\alpha(A, h) = F$;

2) there is a number j such that $\alpha(A, j) = T$ or $\alpha(A, j) = F$.

Let us call a *dated fact* the truth condition of a temporally absolute statement, that is, the (real or unreal) obtaining of a state of affairs at a definite time. In the temporally absolute sense, an utterance now of "Socrates died in Athens", for example, states a real dated fact, because Socrates actually died in Athens at a time before now; and an utterance of the same sentence at a given moment before Socrates' death stated an unreal dated fact because Socrates actually did not die in Athens before that moment. If now it is neither necessary for a sea battle to take place tomorrow nor necessary for it not to take place tomorrow, then the dated fact stated by an utterance now of "There will be a sea battle tomorrow" is now neither determined to be real nor determined to be unreal.

DI9 valuations are intended to be formal proxies for possible worlds. Given any DI9 valuation $\alpha$ and any number j, the set of atomic formulas A of L such that $\alpha(A, j) = T$ stands formally for the set of elementary dated facts that are real at the moment j, or were real until j, or are already determined at j to be real after j in the world $\alpha$; the set of atomic formulas A of L such that $\alpha(A, j) = F$ stands formally for the set of elementary dated facts that are unreal at the moment j, or were unreal until j, or are already determined at j to be unreal after j in the world $\alpha$; and the set of atomic formulas A of L such that $\alpha(A, j) = 0$ stands formally for the set of elementary dated facts that are neither determined at j to be real after j nor determined at j to be unreal after j in the world $\alpha$.

The first condition in Definition 1 stands formally for the thesis that an atomic statement having a truth value at a given moment in a given world will keep this truth value forever after that moment in this world. The rationale behind this thesis is the law of the necessity of the past, which ensures that a real (unreal) past dated fact will remain a real (unreal) past dated fact forever. The second condition stands formally for instances of the law of weak bivalence: in every world, at any moment any atomic statement has or will have a truth value.

*Definition 2*. If $\alpha$ is a DI9 valuation for L, then the *DI9 classical interpretation* for L associated to $\alpha$ is the function $\alpha^*$ from the set of formulas of L to the set of truth-values {T, F} such that:

>1) if A is an atomic formula of L, then $\alpha^*(A) = T$ if and only if there is a j such that $\alpha(A, j) = T$;

>2) If A is ~B, then $\alpha^*(A) = T$ if and only if $\alpha^*(B) = F$;

>3) if A is $B \lor C$, then $\alpha^*(A) = T$ if and only if $\alpha^*(B) = T$ or $\alpha^*(C) = T$.

Conditions 1) and 2) in Definition 1 clearly guarantee the existence and uniqueness of $\alpha^*(A)$, for all formulas A of L. It is worth noting that DI9 classical interpretations are interpretations for L in the sense of classical semantics for classical propositional logic; we will prove later that the converse is also true.

The function $\alpha^*$ is intended to be the formal proxy for the relation between any statement and the only truth value that it has or eventually will have in the world $\alpha$. Proposition 2 below will show that this intention succeeds.

Next, I define the formal counterpart of the concept of a possible continuation of the whole chain of dated facts that were real in the world $\alpha$ until the moment j, including j.

*Definition 3*. Let j be any number. A DI9 valuation $\beta$ for L is a *j-extension* of a DI9 valuation $\alpha$ for L if and only if, for any atomic formula A of L and any number h, if $h \leq j$, then $\beta(A, h) = \alpha(A, h)$.

Finally, I define the notion of a DI9 interpretation for L, which is meant to be a generalization of the notion of a DI9 valuation for L conforming to the logical-semantic theses (1) – (6) above.

*Definition* 4. For any DI9 valuation $\alpha$ for L, the *DI9 interpretation* $I\alpha$ for L is the function from the Cartesian product of the set of formulas of L and the set of real numbers to the set of values {T, F, 0} such that

1) if A is an atomic formula, then $I\alpha(A, j) = \alpha(A, j)$;

2) if A is ~B, then

   2.1) $I\alpha(A, j) = T$ if and only if $I\alpha(B, j) = F$;

   2.2) $I\alpha(A, j) = F$ if and only if $I\alpha(B, j) = T$;

   2.3) $I\alpha(A, j) = 0$ if and only if $I\alpha(B, j) = 0$;

3) if A is $B \vee C$, then

   3.1) $I\alpha(A, j) = T$ if and only if, for any j-extension $\beta$ of $\alpha$, $\beta^*(B) = T$ or $\beta^*(C) = T$;

   3.2) $I\alpha(A, j) = F$ if and only if, for any j-extension $\beta$ of $\alpha$, $\beta^*(B) = F$ and $\beta^*(C) = F$;

   3.3) $I\alpha(A, j) = 0$ if and only if $I\alpha(A, j) \neq T$ and $I\alpha(A, j) \neq F$.

Each DI9 interpretation assigns one and only one value to each formula of L. The crucial step in its recursive definition is the third one, concerning disjunctions. It expresses formally the claim that a disjunction is true in a world at a given moment if and only if in any possible continuation of this world from that moment on at least one disjunct is or will be true; a disjunction is false in a world at a given moment if and only if both disjuncts are or will be false in all possible continuations of this world from that moment on; and a disjunction is neither true nor false in a world at a given moment if and only if (i) at least one disjunct is or will be true in some possible continuation of this world from that moment on, and (ii) both disjuncts are or will be false in some other possible continuation of this world from that moment on.

*Definition* 5. A DI9 interpretation $I\beta$ for L is a *j-extension* of a DI9 interpretation $I\alpha$ for L if and only if, for any formula A of L and any number h, if $h \leq j$, then $I\beta(A, h) = I\alpha(A, h)$.

*Proposition* 1. If a DI9 valuation $\beta$ for L is a *j-extension* of a DI9 valuation $\alpha$ for L, then the DI9 interpretation $I\beta$ for L is a j-extension of the DI9 interpretation $I\alpha$ for L.

*Proof*. Let us assume that a DI9 valuation $\beta$ for L is a j-extension of a DI9 valuation $\alpha$ for L. We prove by induction on the length of A that

(i) $I\beta(A, h) = I\alpha(A, h)$, for any formula A of L and any h such that $h \leq j$.

If A is an atomic formula, then (i) follows trivially from Definitions 3 and 4. If A is a negation, then (i) follows trivially from the inductive hypothesis and Definition 4. Now let A be a disjunction $B \vee C$. It follows immediately from the initial assumption and Definition 3 that the set of j-extensions of $\alpha$ is the set of j-extensions of $\beta$, and it follows immediately from Definition 4 that the value that a DI9 interpretation $I\gamma$ assigns to a disjunction $B \vee C$ is a function of the values that the DI9 classical interpretations associated to the j-extensions of $\gamma$ assign to B and C. Hence, $I\beta(B \vee C, h) = I\alpha(B \vee C, h)$.

It will then be proven that DI9 interpretations fulfill in the whole domain of formulas of L relevant conditions that DI9 valuations fulfill in the restricted domain of atomic formulas of L.

*Proposition 2*. If A is a formula of L, $\alpha$ is a DI9 valuation for L, and W is a truth value, then

(i) $\alpha*(A) = W$ if and only if there is a number j such that $I\alpha(A, j) = W$.

*Proof*. Let A be a formula of L, $\alpha$ be a DI9 valuation for L, and W be a truth value. We prove (i) by induction on the length of A. If A is an atomic formula, then (i) follows trivially from Definitions 2 and 4. If A is a negation, then (i) follows trivially from the inductive hypothesis and Definitions 2 and 4. Now let A be a disjunction $B \vee C$. By Definition 2, $\alpha*(B \vee C) = T$ if and only if $\alpha*(B) = T$ or $\alpha*(C) = T$. By inductive hypothesis, $\alpha*(B) = T$ or $\alpha*(C) = T$ if and only if there is a j such that $I\alpha(B, j) = T$ or $I\alpha(C, j) = T$. By Proposition 1, there is a j such that $I\alpha(B, j) = T$ or $I\alpha(C, j) = T$ if and only if there is a j such that, for any j-extension $\beta$ of $\alpha$, $I\beta(B, j) = T$ or $I\beta(C, j) = T$. By the inductive hypothesis, there is a j such that, for any j-extension β of α, Iβ(B, j) = T or Iβ(C, j) = T if and only if there is a j such that, for any j-extension β of α, β*(B) = T or β*(C) = T. By definition 4, there is a j such that, for any j-extension β of α, β*(B) = T or β*(C) = T if and only if there is a j such that Iα(B ∨ C, j) = T. Hence, α*(B ∨ C) = T if and only if there is a j such that Iα(B ∨ C, j) = T. It can be similarly proven that α*(B ∨ C) = F if and only if there is a j such that Iα(B ∨ C, j) = F.

*Proposition 3*. Let A be a formula of L and α be a DI9 valuation for L. There is a number j such that Iα(A, j) = T or Iα(A, j) = F.

*Proof*. Trivial, by Definition 2 and Proposition 2.

*Proposition 4*. Let A be a formula of L, α be a DI9 valuation for L and W be a truth value. If j and h are numbers such that j ≤ h, then

(i) if Iα(A, j) = W, then Iα(A, h) = W.

*Proof*. Let A be a formula of L, α be a DI9 valuation for L and W be a truth value. Let j and h be numbers such that j ≤ h. We prove (i) by induction on the length of A. If A is an atomic formula, then (i) follows trivially from Definitions 1 and 4. If A is a negation, then (i) follows trivially from the inductive hypothesis and Definition 4. Let A be a disjunction B ∨ C and Iα(B ∨ C) = T. By Definition 4,

(ii) for all j-extensions β of α, β*(B) = T or β*(C) = T.

Let γ be any h-extension of α; by Definition 3, γ is also a j-extension of α, and so, by (ii), γ*(B) = T or γ*(C) = T; therefore, by Definition 4, Iα(A, h) = T. It can be similarly proven that If Iα(A, j) = F, then Iα(A, h) = F.

Proposition 2 shows that the function α* is the formal proxy of the assignment to any statement, either atomic or molecular, of the truth value it has or will have in the world α. Proposition 3 represents formally the fact that in every world at any moment every statement has *or will have* a truth value. Proposition 4 represents formally the fact that a statement that is true (false) in a world at a moment remains true (false) in this world forever after that moment.

Now we prove that DI9 semantics vindicates the definition of truth and falsity of *De Interpretatione* 9.

*Proposition 5*. Let A be a formula of L, α be a DI9 valuation for L, j be any number and W be a truth value; Iα(A, j) = W if and only if, for any j-extension β of α, β*(A) = W.

*Proof*. Let A be a formula of L, α be a DI9 valuation for L, j be any number and W be a truth value. Let Iα(A, j) = W. By Proposition 1, for any j-extension β of α, Iβ(A, j) = W, and so β*(A) = W, by Proposition 2. Now let β*(A) = W, for any j-extension β of α. We prove by induction on the length of A that

(i) Iα(A, j) = W.

If A is a negation or a disjunction, then (i) follows trivially from the inductive hypothesis and Definition 4. Now let A be an atomic formula. We define the functions γ and δ from the Cartesian product of the set of atomic formulas of L and the set of real numbers to the set of values {T, F, 0} as follows: for any atomic formula B of L and any number h such that h ≤ j, γ(B, h) = δ(B, h) = α(B, h); for any number h such that h > j, γ(B, h) = T if and only if α(B, j) = T, and γ(B, h) = F if and only if α(B, j) ≠ T; for any number h such that h > j, δ(B, h) = F if and only if α(B, j) = F, and δ(B, h) = T if and only if α(B, j) ≠ F. It can be easily seen that γ and δ are DI9 valuations for L and j-extensions of α. By definition, for all h > j, γ(A, h) ≠ T if α(A, j) ≠ T, and δ(A, h) ≠ F if α(A, j) ≠ F. Hence, by Definitions 1 and 4, if Iα(A, j) ≠ W, then there is a j-extension β of α such that, for all k, Iβ(A, k) ≠ W.

Proposition 5 represents formally the fact that a statement is true (false) in a given world at a given moment if and only if it is or will be true (false) in all possible continuations of this world from that moment on – in other words, if and only if it is *necessary* at that moment for it to be true (false) in this world.

Finally, we prove that reasonable concepts of logical truth and logical consequence can be defined in DI9 semantics which are coextensive with the concepts of tautology and tautological consequence in classical semantics.

*Definition 6*. A DI9 interpretation Iα for L *satisfies* a formula A of L at a number j if and only if Iα(A, j) = T.

*Definition 7*. A formula A of L is a *DI9 logical consequence* of a set Γ of formulas of L if and only if, for any number j, every DI9 interpretation for L that satisfies all elements of Γ at j also satisfies A at j.

*Definition 8*. A formula A of L is a *DI9 logical truth* if and only if A is a DI9 logical consequence of the empty set.

Definition 8 trivially implies that a formula A of L is a DI9 logical truth if and only if A is satisfied by every DI9 interpretation at every number.

Informally, Definition 6 says that a world satisfies a statement at a given moment if and only if this statement is true in this world at that moment. Definition 7 says that a statement is a DI9 logical consequence of a set of statements if and only if, for any moment, this statement is true at that moment in every world in which all elements of this set are true at that moment. Definition 8 implies that a statement is a DI9 logical truth if and only if it is true in every world at every moment.

*Proposition 6*. If a formula A of L is a tautological consequence of a set Γ of formulas of L in the sense of classical semantics, then A is a DI9 logical consequence of Γ.

*Proof*. Let A be a tautological consequence of Γ, j be any number and $I\alpha$ be any DI9 interpretation that satisfies all elements of Γ at j. By Proposition 1, for all j-extensions β of α and all elements B of Γ, Iβ satisfies B at j and so, by Definition 2, $\beta^*(B) = T$. Since, for all j-extensions β of α, $\beta^*$ is an interpretation in the sense of classical semantics and A is a tautological consequence of Γ, it follows that, for all j-extensions β of α, $\beta^*(A) = T$. By Proposition 5, $I\alpha(A, j) = T$.

*Corollary*. All tautologies are DI9 logical truths.

*Proposition 7*. Let Ic be any interpretation for L in the sense of classical semantics and α be the DI9 valuation for L such that, for all atomic formulas A of L and all numbers j, $\alpha(A, j) = Ic(A)$; we have that

(i) for all formulas B of L and all numbers j, $I\alpha(B, j) = Ic(B)$.

*Proof*. Let Ic be any interpretation for L in the sense of classical semantics and α be the DI9 valuation for L such that, for all atomic formulas A of L and all numbers j,

α(A, j) = Ic(A). We prove (i) by induction on the length of B. If B is an atomic formula, then (i) follows trivially from Definition 4. If B is a negation, then (i) follows trivially from the inductive hypothesis and Definition 4. Let B be a disjunction $C \vee D$. By the inductive hypothesis, for all numbers j, Iα(C, j) = Ic(C) and Iα(D, j) = Ic(D). By Proposition 1, for all numbers j and all j-extensions β of α, Iβ(C, j) = Ic(C) and Iβ(D, j) = Ic(D). By Proposition 2, for all numbers j and all j-extensions β of α, β*(C) = Ic(C) and β*(D) = Ic(D). By Definition 2, for all numbers j and all j-extensions β of α, β*(C $\vee$ D) = Ic(C $\vee$ D). By Proposition 5, for all numbers j, Iα(C $\vee$ D, j) = Ic(C $\vee$ D).

*Proposition 8*. If a formula A of L is a DI9 logical consequence of a set Γ of formulas of L, then A is a tautological consequence of Γ in the sense of classical semantics.

Proof. Let A be a DI9 logical consequence of Γ. Let Ic be any interpretation for L in the sense of classical semantics such that Ic(B) = T, for all elements B of Γ. Let α be the DI9 valuation for L such that, for all atomic formulas C of L and all numbers j, α(C, j) = Ic(C). By Proposition 7, for all formulas D of L and all numbers j, Iα(D, j) = Ic(D). Hence, for all elements B of Γ and all numbers j, Iα(B, j) = Ic(B) = T. Since A is a DI9 logical consequence of Γ, Iα(A, j) = T, for all numbers j, and so Ic(A) = T, by Proposition 7.

*Corollary.* All DI9 logical truths are tautologies.

All DI9 classical interpretations for L are obviously interpretations for L in the sense of classical semantics, and Proposition 7 ensures that every interpretation for L in the sense of classical semantics is a DI9 classical interpretation for L. Therefore, Propositions 6 and 8 imply Proposition 9.

*Proposition 9*. A formula A of L is a DI9 logical consequence of a set of formulas Γ of L if and only if, for all DI9 classical interpretations α* for L, α*(A) = T if, for all elements B of Γ, α*(B) = T.

*Corollary.* A formula A of L is a DI9 logical truth if and only if, for all DI9 classical interpretations α* for L, α*(A) = T.

Informally, Proposition 9 says that a statement is a DI9 logical consequence of a set of statements if and only if it is or will be true in every world in which all elements of this set are or will be true. Its corollary says that a statement is a DI9 logical truth if and only if it is true in every world at some moment. Hence, for a statement to be true in every world at *some* moment is a sufficient (and obviously necessary) condition for it to be true in every world at *every* moment.

*Conclusion*

To sum things up, the scruples of Aristotle's sympathizers in endorsing the traditional reading of *De Interpretatione* 9, as well as the criticisms of his opponents who endorse this reading, are unjustified, for the logical-semantic theses that it finds in the chapter are actually much more conservative than they may appear at first glance.

Concerning the replacement of strong bivalence by weak bivalence, the traditional reading qualifies but does not break the essential link between the concept of statement and the attribute of being true or false. According to it, not every statement is at all times true or false, but every statement that is neither true nor false at a given moment will necessarily be true or false at some other moment. Insofar as a statement is a representation of reality, a kind of representation that can and must be right or wrong, then it is essential for every statement to be at least potentially true or false, and it is also essential for it to be actually true or false at some moment, at the right moment, that is, at the moment in which it is already determined to be real or unreal what it asserts to be real. By force of the weak version of the principle of plenitude, in all cases such a moment will necessarily come.

Concerning the definitions of truth and falsity provided by Aristotle in the *Metaphysics*, the logical semantics of *De Interpretatione* 9 is not disruptive at all. Those definitions disregard distinctions of temporal reference between statements. For that reason, they cannot supply an answer to the issues at stake in the Sea Battle Argument. Since, by definition, the fulfillment or non-fulfillment of the necessary and sufficient condition of truth or falsity of a temporal statement is temporally localized, how is to be temporally localized the truth value of this statement? How truth values should be assigned to statements about the future, before the time when, by definition, the necessary and sufficient conditions of their truth or falsity are to be fulfilled?

In *De Interpretatione* 9, Aristotle supplies a consistent answer to these questions. There he does not contradict but complicates the definitions of the *Metaphysics* in order to make them appropriate to be applied in the domain of temporal statements. He does so by introducing temporally relative concepts of truth and falsity, and defining them through the temporal modalization of the T-scheme and its analogue concerning falsity: S is true (false) at a given moment if and only if it is *necessary* at that moment that S (not-S).

Still, Aristotle preserves a weak version of these schemes: it is *at some moment* true (false) that S if and only if S (not-S). As a consequence, the semantics of *De Interpretatione* 9 is even compatible with a weak variant of extensionalism: the truth value of a molecular statement is *in a way* a function of the truth values of its atomic parts.

In fact, we saw that a crucial step in the definition of the truth conditions of molecular statements in DI9 semantics is the association of truth functions to propositional connectives and the setting of recursive rules by means of which the only truth value that any molecular statement has *or will have* in a given possible world can be assigned solely on the basis of the truth values that its atomic parts have *or will have* in this world. Thus, in the semantics of *De Interpretatione* 9, a disjunction can, to be sure (and to Quine's scandal), be true at a given moment without any of its disjuncts being true *at that moment*, but it can never be true without any of its disjuncts being true *at some moment*.

Finally, the anachronistic formal treatment of classical propositional logic based on DI9 semantics made clear that *De Interpretatione* 9 is no doubt disruptive to classical semantics of classical propositional logic, but it is not disruptive at all to classical propositional logic itself. The adoption of a temporally relative concept of truth applicable to temporally absolute statements and the restriction of the law of strong bivalence are no doubt disruptive on the metalinguistic level, on the level of semantic theory. Nevertheless, they are not disruptive at all on the level of the object language of classical propositional logic itself, as far as they keep untouched the extensions of the classical concepts of propositional logical truth and propositional logical implication. Through different paths, at the end of the story, classical semantics and DI9 semantics eventually agree about what is to be taken as true, and what is to be taken as implied by what, solely in virtue of the meanings of the standard propositional connectives of classical propositional logic – no matter how differently they conceive these meanings to be.